\documentclass{elsart}
\usepackage{amssymb}
\usepackage{amsfonts}
\usepackage{amsmath}

\input{tcilatex}
\begin{document}

\begin{frontmatter}
\title{New Representations of Bertrand Pairs in Euclidean 3-Space}
\author{Y{\i}lmaz Tun{\c{c}}er and Serpil {\"{U}}nal}
\address{Usak University \\ Science and Arts Faculty \\ Department of Mathematics \\ Corresponding Author: yilmaz.tuncer@usak.edu.tr}
\begin{abstract}
In this work, we studied the properties of the spherical indicatrices of a Bertrand curve and its mate curve and presented some characteristic properties in the cases that Bertrand curve and its mate curve are slant helices, spherical indicatrices are slant helices and we also researched that whether the spherical indicatrices made new curve pairs in the means of Mannheim, involte-evolute and Bertrand pairs. Further more, we investigated the relations between the spherical images and introduced new representations of spherical indicatrices.
\end{abstract}
\begin{keyword}
Bertrand curve, Bertrand mate curve, Spherical indicatrix, Slant helix.\\ \textit{AMS Subject Classification 2000}: 53A04
\end{keyword}
\end{frontmatter}

\section{Introduction}

\bigskip It is well known that the specific curve pairs are the most popular
subjects in curve and surface theory and Bertrand curve is one of them. A
Bertrand curve is a curve in Euclidean 3-space whose principal normal is the
principal normal of another curve. We can see in most textbooks, a
characteristic property of the Bertrand curve which asserts the existence of
a linear relation between curvature and torsion. The characteristic
property-linear relation-is deduced as an application of the Frenet-Serret
formulae.

In this study, the spherical indicatrices of a Bertrand curve and its mate
curve are given. In order to make a Bertrand curve and its mate curve slant
helix, the feature that spherical indicatrices curve's need to have are
examined. It is seen that the curves of the spherical indicatrices do not
form any specific curve pairs.

Let $\gamma :I\longrightarrow IR^{3}$ be a curve with $\gamma ^{\prime
}\left( s\right) \neq 0$, where $\gamma ^{\prime }\left( s\right) =d$ $%
\gamma \left( s\right) /ds$. The arc-lenght $s$ of a curve $\gamma \left(
s\right) $ is determined such that $\left\Vert \gamma ^{\prime }\left(
s\right) \right\Vert =1.$ This shows us that $T\left( s\right) =\gamma
^{\prime }\left( s\right) $ and $T\left( s\right) $ is called a tangent
vector of $\gamma $ at $\gamma (s)$. We define the curvature of $\gamma $ by 
$\kappa \left( s\right) =\left\Vert \gamma ^{\prime \prime }\left( s\right)
\right\Vert $. If $\kappa \left( s\right) \neq 0,$ the unit principal normal
vector $N\left( s\right) $ of the curve at $\gamma (s)$ is given by $\gamma
^{\prime \prime }\left( s\right) =\kappa \left( s\right) N\left( s\right) $.
The unit vector $B\left( s\right) =T\left( s\right) \Lambda N\left( s\right) 
$ is called the unit binormal vector of $\gamma $ at $\gamma (s)$. Then we
have the Frenet-Serret formulae 
\begin{equation*}
T^{\prime }=\kappa N\text{ , \ }N^{\prime }=-\kappa T+\tau B\text{ , \ }%
B^{\prime }=-\tau N
\end{equation*}%
where $\tau \left( s\right) $ is the torsion of $\gamma $ at $\gamma (s)$ 
\cite{5}.

The curve $\gamma $ with nonzero curvature is called a Bertrand curve if
there exists a curve $\widetilde{\gamma }:I\longrightarrow IR^{3}$ such that
the principal normal vectors of $\gamma $ and $\widetilde{\gamma }$ are
linearly dependent at the corresponding points for each $s\in I\subset IR.$
In this case, $\widetilde{\gamma }$ is called a Bertrand mate of $\gamma $
and there exists a relationship between the position vectors as%
\begin{equation*}
\widetilde{\gamma }(s^{\ast })=\gamma (s)+\lambda N(s)
\end{equation*}%
and we can write $N=\epsilon \widetilde{N}$ where $\epsilon =\pm 1$ and $%
\lambda $ is the distance between the curves $\gamma $ and $\widetilde{%
\gamma }$ at the corresponding points for each $s.$ The pair of ($\gamma $, $%
\widetilde{\gamma }$) is called a Bertrand pair. $\lambda $ is a constant
for Bertrand pairs\cite{3,6}.

The curve $\widetilde{\gamma }$ is called involute of $\gamma $ if the
tangent vectors are orthogonal at the corresponding points for each $s\in
I\subset IR.$ In this case, $\gamma $ is called evolute of the curve $%
\widetilde{\gamma }$ and there exists a relationship between the position
vectors as%
\begin{equation*}
\widetilde{\gamma }(s^{\ast })=\gamma (s)+\lambda T(s)
\end{equation*}%
where $\lambda $ is the distance between the curves $\gamma $ and $%
\widetilde{\gamma }$ at the corresponding points for each $s.$The pair of ($%
\gamma $, $\widetilde{\gamma }$) is called a involute-evolute pair. $\lambda 
$ is not a constant for involute-evolute pairs\cite{5}.

The curve $\widetilde{\gamma }$ is called Mannheim curve if there exists a
corresponding relationship between the space curves $\gamma $ and $%
\widetilde{\gamma }$ such that the principal normal lines of $\widetilde{%
\gamma }$ coincide with the binormal lines of $\gamma $ at the corresponding
points of the curves. In this case, $\gamma $ is called as a Mannheim
partner curve of $\widetilde{\gamma }$ and there exists a relationship
between the position vectors as%
\begin{equation*}
\widetilde{\gamma }(s^{\ast })=\gamma (s)+\lambda B(s)
\end{equation*}%
where $\lambda $ is the distance between the curves $\gamma $ and $%
\widetilde{\gamma }$ at the corresponding points for each $s.$ The pair of ($%
\gamma $, $\widetilde{\gamma }$) is called a Mannheim pair. $\lambda $ is a
constant for Mannheim pairs\cite{1,2}.

On the other hand, Izumiya and Takeuchi have introduced the concept of slant
helix by saying that the normal lines make a constant angle with a fixed
straight line. They characterize a slant helix if and only if the geodesic
curvature of the principal image of the principal normal indicatrix%
\begin{equation*}
\Gamma =\frac{\kappa ^{2}}{\left( \kappa ^{2}+\tau ^{2}\right) ^{3/2}}\left( 
\frac{\tau }{\kappa }\right) ^{\prime }
\end{equation*}%
is a constant function\cite{4,6} .

In this study, we denote that $T$, $N$, $B,$ $\kappa $, $\tau $ and $%
\widetilde{T}$, $\widetilde{N}$, $\widetilde{B},$ $\widetilde{\kappa }$, $%
\widetilde{\tau }$ are the Frenet equipments of $\gamma $ and $\widetilde{%
\gamma },$ respectively. Tangent, principal normal and binormal vectors are
described for the spherical curves which are called tangent, principal
normal and binormal indicatrices for both the curves $\gamma $ and $%
\widetilde{\gamma },$ respectively. Throughout this study, the Bertrand
curve and its mate curve are accepted as non helix and non planar.

\section{Spherical indicatrices of a Bertrand Curve}

\bigskip In this section, we introduced the spherical indicatrices of a
Bertrand curve according to its mate curve in Euclidean 3-space and gave
considerable results by using the properties of the curves. Let $\gamma $ be
a Bertrand curve with a Bertrand mate curve $\widetilde{\gamma ,}$ then 
\begin{equation}
\gamma (s)=\widetilde{\gamma }(s^{\ast })-\lambda \epsilon \widetilde{N}
\label{1}
\end{equation}%
where $\lambda $ is 
\begin{equation}
\lambda =\frac{-\epsilon \widetilde{g}}{\widetilde{\kappa }\left( \widetilde{%
g}-\widetilde{f}\right) }  \label{2}
\end{equation}%
and 
\begin{equation*}
\widetilde{f}=\frac{\widetilde{\tau }}{\widetilde{\kappa }},\text{ \ \ }%
\widetilde{g}=\frac{\widetilde{\tau }^{\prime }}{\widetilde{\kappa }^{\prime
}}.
\end{equation*}

\begin{thm}
Let $\gamma $ be a Bertrand curve, then we have Frenet formula:%
\begin{equation*}
T^{\prime }=\kappa N,\text{ \ \ }N^{\prime }=-\kappa T+\tau B,\text{ \ \ }%
B^{\prime }=-\tau N
\end{equation*}%
where%
\begin{equation}
T=\frac{-1}{\sqrt{1+\widetilde{g}^{2}}}\left\{ \widetilde{T}-\widetilde{g}%
\widetilde{B}\right\} ,\text{ \ \ }N=\epsilon \widetilde{N},\text{ \ \ }B=%
\frac{-\epsilon }{\sqrt{1+\widetilde{g}^{2}}}\left\{ \widetilde{g}\widetilde{%
T}+\widetilde{B}\right\}  \label{3}
\end{equation}%
and%
\begin{equation}
s=\dint \frac{\widetilde{f}\sqrt{1+\widetilde{g}^{2}}}{\widetilde{g}-%
\widetilde{f}}ds^{\ast },\text{ \ \ }\kappa =\frac{-\epsilon \widetilde{%
\kappa }\left( 1+\widetilde{f}\widetilde{g}\right) \left( \widetilde{g}-%
\widetilde{f}\right) }{\widetilde{f}\left( 1+\widetilde{g}^{2}\right) },%
\text{ \ \ }\tau =\frac{\widetilde{\kappa }\left( \widetilde{g}-\widetilde{f}%
\right) ^{2}}{\widetilde{f}\left( 1+\widetilde{g}^{2}\right) }  \label{4}
\end{equation}%
The geodesic curvature of the principal image of the principal normal
indicatrix of $\gamma $ is%
\begin{equation}
\Gamma =\frac{-\widetilde{\kappa }^{\prime }\left( \widetilde{g}-\widetilde{f%
}\right) }{\widetilde{\kappa }^{2}\left( 1+\widetilde{f}^{2}\right) ^{3/2}}
\label{5}
\end{equation}
\end{thm}

\begin{thm}
\label{th2}A Bertrand curve is a$\ $slant helix if and only if its mate
curve is a slant helix.
\end{thm}

\begin{pf}
Let $\Gamma $ and $\widetilde{\Gamma }$ be the geodesic curvatures of the
principal normal curves of$\ \gamma $ and $\widetilde{\gamma },$
respectively. Then $\widetilde{\Gamma }$ is 
\begin{equation*}
\widetilde{\Gamma }=\frac{\widetilde{\kappa }^{\prime }\left( \widetilde{g}-%
\widetilde{f}\right) }{\widetilde{\kappa }^{2}\left( 1+\widetilde{f}%
^{2}\right) ^{3/2}}
\end{equation*}%
and from (\ref{5}) we get 
\begin{equation*}
\Gamma =-\widetilde{\Gamma }
\end{equation*}%
So, it is easy to see that $\gamma $ is a slant helix if and only if $%
\widetilde{\gamma }$ is a slant helix.
\end{pf}

Let $p_{1}=\frac{1}{\sqrt{1+\widetilde{g}^{2}}}\ $and$\ p_{2}=\frac{%
\widetilde{g}}{\sqrt{1+\widetilde{g}^{2}}}.$ If we differentiate the
equation (\ref{1}) with respect to $s^{\ast },$ we get%
\begin{equation*}
\gamma ^{\prime }(s)\frac{ds}{ds^{\ast }}=\frac{-\widetilde{f}}{\widetilde{g}%
-\widetilde{f}}\widetilde{T}+\frac{\widetilde{f}\widetilde{g}}{\widetilde{g}-%
\widetilde{f}}\widetilde{B}.
\end{equation*}%
Since%
\begin{equation*}
\frac{ds^{\ast }}{ds}=\frac{\widetilde{g}-\widetilde{f}}{\widetilde{f}\sqrt{%
1+\widetilde{g}^{2}}},
\end{equation*}%
we can write%
\begin{equation*}
\gamma ^{\prime }(s)=p_{1}\widetilde{T}+p_{2}\widetilde{B}
\end{equation*}%
again, if we differentiate the last equation with respect to $s^{\ast }$ and
put $\lambda $\ from equation (\ref{2}), we obtain%
\begin{equation*}
\gamma ^{\prime \prime }(s)=-p_{1}^{\prime }\frac{ds^{\ast }}{ds}\widetilde{T%
}-\left\{ p_{1}\widetilde{\kappa }+p_{2}\widetilde{\tau }\right\} \frac{%
ds^{\ast }}{ds}\widetilde{N}+p_{2}^{\prime }\frac{ds^{\ast }}{ds}\widetilde{B%
}.
\end{equation*}%
Since $\gamma ^{\prime \prime }(s)$ has no components on $\widetilde{T}$ and 
$\widetilde{B}$, $p_{1}^{\prime }=0$ and $p_{2}^{\prime }=0.$ Hence, $p_{1}$
and $p_{2}$ are constants. Thus, we prove the following theorem.

\begin{thm}
\label{th3}Let $\gamma $ and $\widetilde{\gamma }$ be regular curves in $%
E^{3}.$If $\gamma $ is a Bertrand curve with its mate curve $\widetilde{%
\gamma },$ $\widetilde{g}$ is a constant.
\end{thm}

Let $\gamma $ be\ a unit speed regular curve in Euclidean 3-space with
Frenet vectors $T$, $N$ and $B$. The unit tangent vectors along the curve $%
\gamma (s)$ generates a curve $\gamma _{t}=T$ on the sphere of radius 1
about the origin. The curve $\gamma _{t}$ is called the spherical indicatrix
of $T$ or more commonly, $\gamma _{t}$ is called the tangent indicatrix of
the curve $\gamma $. If $\gamma =\gamma (s)$ is a natural representations of
the curve $\gamma $, $\gamma _{t}(s)=T(s)$ will be a representation of $%
\gamma _{t}$. Similarly, we can consider the principal normal indicatrix $%
\gamma _{n}=N(s)$ and the binormal indicatrix $\gamma _{b}=B(s).$

The tangent indicatrix of a Bertrand curve is%
\begin{equation*}
\gamma _{t}=\frac{-1}{\sqrt{1+\widetilde{g}^{2}}}\left\{ \widetilde{T}-%
\widetilde{g}\widetilde{B}\right\} .
\end{equation*}

\begin{thm}
If the Frenet frame of the tangent indicatrix $\gamma _{t}=T$ of the
Bertrand curve is $\{T_{t},N_{t},B_{t}\}$, we have Frenet formula:%
\begin{equation*}
T_{t}^{\prime }=\kappa _{t}N_{t},\text{ \ \ }N_{t}^{\prime }=-\kappa
_{t}T_{t}+\tau _{t}B_{t},\text{ \ \ }B_{t}^{\prime }=-\tau _{t}N_{t}
\end{equation*}%
where%
\begin{equation}
T_{t}=-\widetilde{N},\text{ \ \ }N_{t}=\frac{1}{\sqrt{1+\widetilde{f}^{2}}}%
\left\{ \widetilde{T}-\widetilde{f}\widetilde{B}\right\} ,\text{ \ \ }B_{t}=%
\frac{1}{\sqrt{1+\widetilde{f}^{2}}}\left\{ \widetilde{f}\widetilde{T}+%
\widetilde{B}\right\}  \label{6}
\end{equation}%
and%
\begin{equation}
s_{t}=-\dint \frac{\widetilde{\kappa }\left( \widetilde{g}-\widetilde{f}%
\right) ^{2}}{\widetilde{f}\left( 1+\widetilde{g}\right) ^{2}}ds,\text{ \ \ }%
\kappa _{t}=\frac{\sqrt{1+\widetilde{g}^{2}}\sqrt{1+\widetilde{f}^{2}}}{%
\left( \widetilde{f}-\widetilde{g}\right) },\ \ \ \tau _{t}=\frac{-%
\widetilde{\kappa }^{\prime }\sqrt{1+\widetilde{g}^{2}}}{\widetilde{\kappa }%
^{2}\left( 1+\widetilde{f}^{2}\right) }  \label{7}
\end{equation}%
and $s_{t}$ is a natural representation of the tangent indicatrix of the
curve $\gamma $. $\kappa _{t}$ and $\tau _{t}$ are the curvature and torsion
of $\gamma _{t}$, respectively. The geodesic curvature of the principal
image of the principal normal indicatrix of $\gamma _{t}$ is%
\begin{equation}
\Gamma _{t}=\frac{-\widetilde{\kappa }^{3}\left( 1+\widetilde{f}^{2}\right)
^{3/2}\left( \widetilde{g}-\widetilde{f}\right) ^{2}\left\{ \widetilde{%
\kappa }^{\prime \prime }\widetilde{\kappa }\left( 1+\widetilde{f}%
^{2}\right) -3\widetilde{\kappa }^{\prime 2}\left( 1+\widetilde{f}\widetilde{%
g}\right) \right\} }{\sqrt{1+\widetilde{g}^{2}}\left( \widetilde{\kappa }%
^{4}\left( 1+\widetilde{f}^{2}\right) ^{3}+\widetilde{\kappa }^{\prime
2}\left( \widetilde{f}-\widetilde{g}\right) ^{2}\right) ^{3/2}}\frac{%
ds^{\ast }}{ds_{t}}.  \label{8}
\end{equation}
\end{thm}

The arclengths $s_{t}$ and $s^{\ast }$\ satisfy the equation%
\begin{equation*}
s_{t}=\dint \frac{\widetilde{\kappa }\left( \widetilde{f}-\widetilde{g}%
\right) }{\sqrt{1+\widetilde{g}^{2}}}ds^{\ast }.
\end{equation*}

We can give the following theorem to characterize a Bertrand curve as a
slant helix.

\begin{thm}
\label{th6}Let $\left( \gamma ,\widetilde{\gamma }\right) $ be a non-helical
and non-planar Bertrand pair parameterized by the arclengths $s$ and $%
s^{\ast }.$Thus, the followings are true.

i. Bertrand curve is a slant helix if and only if the tangent indicatrix of
the Bertrand curve is a spherical helix.

ii. Bertrand mate curve is a slant helix if and only if the tangent
indicatrix of the Bertrand curve is a spherical helix.
\end{thm}

Thus, we can state the following theorem to characterize the spherical image
of tangent indicatrix of a Bertrand curve as a helix.

\begin{thm}
\label{th8}Let $\left( \gamma ,\widetilde{\gamma }\right) $ be a non-helical
and non-planar Bertrand pair$.$Thus, the tangent indicatrix of the Bertrand
curve is a spherical helix if and only if 
\begin{equation*}
\widetilde{\kappa }^{\prime \prime }\widetilde{\kappa }\widetilde{f}^{2}-3%
\widetilde{\kappa }^{\prime 2}\widetilde{g}\widetilde{f}+\widetilde{\kappa }%
^{\prime \prime }\widetilde{\kappa }-3\widetilde{\kappa }^{\prime 2}=0
\end{equation*}%
is satisfied.
\end{thm}

The principal normal indicatrix of the Bertrand curve is%
\begin{equation*}
\gamma _{n}=\epsilon \widetilde{N}.
\end{equation*}

\begin{thm}
If the Frenet frame of the principal normal indicatrix $\gamma _{n}=N$ of
the Bertrand curve is $\{T_{n},N_{n},B_{n}\}$, we have Frenet formula:%
\begin{equation*}
T_{n}^{\prime }=\kappa _{n}N_{n},\text{ \ \ }N_{n}^{\prime }=-\kappa
_{n}T_{n}+\tau _{n}B_{n},\text{ \ \ }B_{n}^{\prime }=-\tau _{n}N_{n}
\end{equation*}%
where%
\begin{eqnarray}
T_{n} &=&\frac{-\epsilon }{\sqrt{1+\widetilde{f}^{2}}}\left\{ \widetilde{T}-%
\widetilde{f}\widetilde{B}\right\}  \notag \\
N_{n} &=&\frac{\epsilon }{\rho \sqrt{1+\widetilde{f}^{2}}}\left\{ \widetilde{%
f}\widetilde{\kappa }^{\prime }\left( \widetilde{g}-\widetilde{f}\right) 
\widetilde{T}-\widetilde{\kappa }^{2}\left( 1+\widetilde{f}^{2}\right) ^{2}%
\widetilde{N}+\widetilde{\kappa }^{\prime }\left( \widetilde{g}-\widetilde{f}%
\right) \widetilde{B}\right\}  \label{9} \\
B_{n} &=&\frac{1}{\rho }\left\{ \widetilde{\kappa }^{2}\widetilde{f}\left( 1+%
\widetilde{f}^{2}\right) \widetilde{T}+\widetilde{\kappa }^{\prime }\left( 
\widetilde{g}-\widetilde{f}\right) \widetilde{N}+\widetilde{\kappa }%
^{2}\left( 1+\widetilde{f}^{2}\right) \widetilde{B}\right\}  \notag
\end{eqnarray}%
and%
\begin{gather*}
s_{n}=\dint \frac{\widetilde{\kappa }\left( \widetilde{g}-\widetilde{f}%
\right) \sqrt{1+\widetilde{f}^{2}}}{\widetilde{f}\sqrt{1+\widetilde{g}^{2}}}%
ds,\text{ \ \ }\kappa _{n}=\frac{\rho }{\widetilde{\kappa }^{2}\left( 1+%
\widetilde{f}^{2}\right) ^{3/2}} \\
\tau _{n}=\frac{-\epsilon \left( \widetilde{g}-\widetilde{f}\right) }{\rho
^{2}}\left\{ \left( 3\widetilde{\kappa }^{\prime 2}-\widetilde{\kappa }%
\widetilde{\kappa }^{\prime \prime }\right) \left( 1+\widetilde{f}%
^{2}\right) +3\widetilde{f}\widetilde{\kappa }^{\prime 2}\left( \widetilde{g}%
-\widetilde{f}\right) \right\}
\end{gather*}%
where%
\begin{equation*}
\rho =\sqrt{\widetilde{\kappa }^{\prime 2}\left( \widetilde{g}-\widetilde{f}%
\right) ^{2}+\widetilde{\kappa }^{4}\left( 1+\widetilde{f}^{2}\right) ^{3}}
\end{equation*}%
and $s_{n}$ is a natural representation of the principal normal indicatrix
of the curve $\gamma $. $\kappa _{n}$ and $\tau _{n}$ are the curvature and
torsion of $\gamma _{n}$, respectively.
\end{thm}

The arclengths $s_{n}$ and $s^{\ast }$\ satisfy the equation%
\begin{equation*}
s_{n}=\dint \widetilde{\kappa }\sqrt{1+\widetilde{f}^{2}}ds^{\ast }.
\end{equation*}

We can give the following theorem to characterize the spherical image of the
principal normal indicatrix of a Bertrand curve as a helix.

\begin{thm}
\label{th11}Let $\left( \gamma ,\widetilde{\gamma }\right) $ be a
non-helical and non-planar Bertrand pair. Thus, the principal normal
indicatrix of the Bertrand curve is planar if and only if%
\begin{equation*}
\widetilde{\kappa }\widetilde{\kappa }^{\prime \prime }\widetilde{f}^{2}-3%
\widetilde{\kappa }^{\prime 2}\widetilde{g}\widetilde{f}-\left( 3\widetilde{%
\kappa }^{\prime 2}-\widetilde{\kappa }\widetilde{\kappa }^{\prime \prime
}\right) =0
\end{equation*}%
is satisfied.
\end{thm}

The binormal indicatrix of the Bertrand curve is%
\begin{equation*}
\gamma _{b}=\frac{-\epsilon }{\sqrt{1+\widetilde{g}^{2}}}\left\{ \widetilde{g%
}\widetilde{T}+\widetilde{B}\right\} .
\end{equation*}

\begin{thm}
\label{elf}If the Frenet frame of the binormal indicatrix $\gamma _{b}=B$ of
the Bertrand curve is $\{T_{b},N_{b},B_{b}\}$, we have Frenet formula:%
\begin{equation*}
T_{b}^{\prime }=\kappa _{b}N_{b},\text{ \ \ }N_{b}^{\prime }=-\kappa
_{b}T_{b}+\tau _{b}B_{b},\text{ \ \ }B_{b}^{\prime }=-\tau _{b}N_{b}
\end{equation*}%
where%
\begin{equation}
T_{b}=\epsilon \widetilde{N},\text{ \ \ }N_{b}=\frac{-\epsilon }{\sqrt{1+%
\widetilde{f}^{2}}}\left\{ \widetilde{T}-\widetilde{f}\widetilde{B}\right\} ,%
\text{ \ \ }B_{b}=\frac{1}{\sqrt{1+\widetilde{f}^{2}}}\left\{ \widetilde{f}%
\widetilde{T}+\widetilde{B}\right\}  \label{10}
\end{equation}%
and%
\begin{equation}
s_{b}=-\dint \frac{\widetilde{\kappa }\left( \widetilde{f}-\widetilde{g}%
\right) ^{2}}{\widetilde{f}\left( 1+\widetilde{g}^{2}\right) }ds,\text{ \ \ }%
\kappa _{b}=\frac{\sqrt{1+\widetilde{f}^{2}}\sqrt{1+\widetilde{g}^{2}}}{%
\widetilde{f}-\widetilde{g}},\ \ \ \tau _{b}=\frac{-\epsilon \widetilde{%
\kappa }^{\prime }\sqrt{1+\widetilde{g}^{2}}}{\widetilde{\kappa }^{2}\left(
1+\widetilde{f}^{2}\right) }  \label{11}
\end{equation}%
and $s_{b}$ is a natural representation of the binormal indicatrix of the
curve $\gamma $ and it is equal to the total torsion of the curve $\gamma $. 
$\kappa _{b}$ and $\tau _{b}$ are the curvature and torsion of $\gamma _{b}$%
, respectively. The geodesic curvature of the principal image of the
principal normal indicatrix of $\gamma _{b}$ is%
\begin{equation}
\Gamma _{b}=\frac{-\widetilde{\kappa }^{3}\left( 1+\widetilde{f}^{2}\right)
^{3/2}\left( \widetilde{g}-\widetilde{f}\right) ^{2}\left\{ \widetilde{%
\kappa }^{\prime \prime }\widetilde{\kappa }\left( 1+\widetilde{f}%
^{2}\right) -3\widetilde{\kappa }^{\prime 2}\left( 1+\widetilde{f}\widetilde{%
g}\right) \right\} }{\sqrt{1+\widetilde{g}^{2}}\left( \widetilde{\kappa }%
^{4}\left( 1+\widetilde{f}^{2}\right) ^{3}+\widetilde{\kappa }^{\prime
2}\left( \widetilde{f}-\widetilde{g}\right) ^{2}\right) ^{3/2}}\frac{%
ds^{\ast }}{ds_{b}}.  \label{12}
\end{equation}
\end{thm}

From theorem \ref{th2} and from (\ref{7}), (\ref{11}), we conclude that 
\begin{equation*}
\Gamma =-\widetilde{\Gamma }=\frac{\tau _{t}}{\kappa _{t}}=\frac{\tau _{b}}{%
\kappa _{b}}
\end{equation*}%
and $\Gamma _{t}=\Gamma _{b}.$Thus, we can give the following corollary.

\begin{cor}
The spherical images of the tangent and binormal indicatrices of a Bertrand
curve are the curves with same curvature and same torsion.
\end{cor}

\begin{cor}
\label{cr14}The arclengths $s_{b}$ and $s^{\ast }$\ satisfy the following
equation.%
\begin{equation*}
s_{b}=\dint \frac{\widetilde{\kappa }\left( \widetilde{f}-\widetilde{g}%
\right) }{\sqrt{1+\widetilde{g}^{2}}}ds^{\ast }
\end{equation*}
\end{cor}

Since $\lambda $ is a constant, from (\ref{2}) we can write 
\begin{equation*}
\widetilde{\tau }^{\prime }=-\epsilon \lambda \widetilde{\kappa }^{2}%
\widetilde{f}^{\prime }
\end{equation*}%
and since $p_{2}$\ is a constant , we can write 
\begin{equation*}
\widetilde{\tau }^{\prime }=c_{1}\widetilde{\kappa }^{\prime }\sqrt{1+%
\widetilde{g}^{2}}.
\end{equation*}%
Thus, we obtain%
\begin{equation*}
\frac{\widetilde{\kappa }^{2}\widetilde{f}^{\prime }}{\widetilde{\kappa }%
^{\prime }\sqrt{1+\widetilde{g}^{2}}}=\frac{-\epsilon c_{1}}{\lambda }
\end{equation*}%
and we can also say that%
\begin{equation*}
\frac{\widetilde{\kappa }^{2}\widetilde{f}^{\prime }}{\widetilde{\kappa }%
^{\prime }\sqrt{1+\widetilde{g}^{2}}}
\end{equation*}%
is a constant. From corollary\ \ref{cr14},we obtain%
\begin{equation*}
\frac{ds_{b}}{ds^{\ast }}=\frac{-\epsilon c_{1}}{\lambda }
\end{equation*}%
and 
\begin{equation*}
s_{b}=\frac{-\epsilon c_{1}}{\lambda }s^{\ast }+c_{2}
\end{equation*}%
where $c_{1}$ and $c_{2}$ are constants, too.

We can give the following theorem to characterize a Bertrand mate curve as a
slant helix.

\begin{thm}
\label{teo15}Let $\left( \gamma ,\widetilde{\gamma }\right) $ be a
non-helical and non-planar Bertrand pair$.$Thus, the Bertrand mate curve is
a slant helix if and only if the binormal indicatrix of the Bertrand curve
is a spherical helix.
\end{thm}

\begin{pf}
Let $\Gamma $ be the geodesic curvature of the principal image of the
principal normal indicatrix of the Bertrand curve. From (\ref{5}) and (\ref%
{11}), we get%
\begin{equation*}
\Gamma =\frac{\tau _{b}}{\kappa _{b}}.
\end{equation*}%
Thus, it is easy to see that the Bertrand mate curve is a slant helix if and
only if the binormal indicatrix of the Bertrand curve is a spherical helix.
\end{pf}

We can give the following corollary by using theorem \ref{th2}.

\begin{cor}
A non-helical and non-planar Bertrand curve is a slant helix if and only if\
the binormal indicatrix of the Bertrand curve is a spherical helix.
\end{cor}

Thus, we can give the following theorem to characterize the spherical image
of the binormal indicatrix of a Bertrand curve as a helix.

\begin{thm}
\label{th17}Let $\left( \gamma ,\widetilde{\gamma }\right) $ be a
non-helical and non-planar Bertrand pair$.$Thus, the binormal indicatrix of
the Bertrand curve is a spherical helix if and only if%
\begin{equation*}
\widetilde{\kappa }^{\prime \prime }\widetilde{\kappa }\widetilde{f}^{2}-3%
\widetilde{\kappa }^{\prime 2}\widetilde{g}\widetilde{f}+\widetilde{\kappa }%
^{\prime \prime }\widetilde{\kappa }-3\widetilde{\kappa }^{\prime 2}=0
\end{equation*}%
is satisfied.
\end{thm}

\begin{cor}
\label{cr18}Let $\left( \gamma ,\widetilde{\gamma }\right) $ be a
non-helical and non-planar Bertrand pair, thus, the followings are
equivalent.

\textbf{i.} The tangent indicatrix of the Bertrand curve is a spherical
helix.

\textbf{ii.} The principal normal indicatrix of the Bertrand curve is planar.

\textbf{iii.} The binormal indicatrix of the Bertrand curve is a spherical
helix.
\end{cor}

The equivalence in Corollary \ref{cr18} can easily prove such as $%
i\Rightarrow ii,ii\Rightarrow iii$ and $iii\Rightarrow i$. We can give the
following corollary by using (\ref{7}), (\ref{9}) and (\ref{10})

\begin{cor}
The equations between the Frenet vectors of the spherical indicatrices of a
non-helical and non-planar Bertrand curve in Euclidean 3-space are $\ \ $%
\begin{equation*}
T_{t}=-\epsilon T_{b},\text{ }T_{n}=-\epsilon N_{t}=N_{b}\text{\ and }%
B_{t}=B_{b}.
\end{equation*}
\end{cor}

\section{Spherical indicatrices of Bertrand Mate Curve}

In this section, we introduced the spherical indicatrices of a Bertrand mate
curve according to the Bertrand curve in Euclidean 3-space and gave
considerable results by using the properties of the curves, similar to the
previous section. Let $\gamma $ be a Bertrand curve with a Bertrand mate
curve $\widetilde{\gamma }$. Then we can write the following equation by
using (\ref{1}) 
\begin{equation}
\widetilde{\gamma }(s^{\ast })=\gamma (s)+\lambda N(s)  \label{13}
\end{equation}%
where $\lambda $ is 
\begin{equation}
\lambda =\frac{g}{\kappa \left( g-f\right) }  \label{14}
\end{equation}%
and%
\begin{equation*}
f=\frac{\tau }{\kappa },\text{ \ \ }g=\frac{\tau ^{\prime }}{\kappa ^{\prime
}}.
\end{equation*}%
We used the same $\lambda $ in (\ref{1}) to see the relation between (\ref{2}%
) and (\ref{14}) such that%
\begin{equation*}
\frac{g}{\kappa \left( g-f\right) }=\frac{-\epsilon \widetilde{g}}{%
\widetilde{\kappa }\left( \widetilde{g}-\widetilde{f}\right) }.
\end{equation*}%
So, we can state that the following equation is satisfied for any
non-helical and non-planar Bertrand curve and its mate curve which is
non-helical and non-planar too. 
\begin{equation*}
\left( \widetilde{\kappa }+\epsilon \kappa \right) g\widetilde{g}-\epsilon f%
\widetilde{g}\kappa -\widetilde{f}g\widetilde{\kappa }=0
\end{equation*}

\begin{thm}
Let $\widetilde{\gamma }$ be a Bertrand mate curve, then we have Frenet
formula:%
\begin{equation*}
\widetilde{T}^{\prime }=\widetilde{\kappa }\widetilde{N},\text{ \ \ }%
\widetilde{N}^{\prime }=-\widetilde{\kappa }\widetilde{T}+\widetilde{\tau }%
\widetilde{B},\text{ \ \ }\widetilde{B}^{\prime }=-\widetilde{\tau }%
\widetilde{N}
\end{equation*}%
where%
\begin{equation*}
\widetilde{T}=-\frac{1}{\sqrt{1+g^{2}}}\left\{ T-gB\right\} ,\text{ \ }\ 
\widetilde{N}=\epsilon N,\ \ \ \widetilde{B}=-\frac{\epsilon }{\sqrt{1+g^{2}}%
}\left\{ gT+B\right\}
\end{equation*}%
and%
\begin{equation*}
ds^{\ast }=\dint \frac{f\sqrt{1+g^{2}}}{\left( g-f\right) }ds,\ \ \ 
\widetilde{\kappa }=\frac{-\epsilon \kappa \left( g-f\right) (1+fg)}{f\left(
1+g^{2}\right) },\text{ \ \ }\widetilde{\tau }=\frac{\kappa \left(
g-f\right) ^{2}}{f\left( 1+g^{2}\right) }.
\end{equation*}%
The geodesic curvature of the the principal image of the principal normal
indicatrix of Bertrand mate curve is%
\begin{equation}
\widetilde{\Gamma }=\frac{\kappa ^{\prime }f\left( 1+g^{2}\right) ^{2}}{%
-\kappa ^{2}\left( (1+fg)^{2}+\left( g-f\right) ^{2}\right) ^{3/2}}\frac{ds}{%
ds^{\ast }}.  \label{15}
\end{equation}
\end{thm}

Let $q_{1}=\frac{1}{\sqrt{1+g^{2}}},$ and $q_{2}=\frac{g}{\sqrt{1+g^{2}}}$.\
If we differentiate the equation (\ref{13}) with respect to $s$, we get%
\begin{equation*}
\widetilde{\gamma }^{\prime }(s^{\ast })\frac{ds^{\ast }}{ds}=\frac{-f}{%
\left( g-f\right) }T+\frac{fg}{\left( g-f\right) }B.
\end{equation*}%
Since%
\begin{equation*}
\frac{ds}{ds^{\ast }}=\frac{\left( g-f\right) }{f\sqrt{1+g^{2}}},
\end{equation*}%
we can write%
\begin{equation*}
\widetilde{\gamma }^{\prime }(s^{\ast })=-q_{1}T+q_{2}B
\end{equation*}%
again, if we differentiate the last equation with respect to $s,$\ we obtain%
\begin{equation*}
\widetilde{\gamma }^{\prime \prime }(s^{\ast })=-q_{1}^{\prime }\frac{ds}{%
ds^{\ast }}T-\kappa \left\{ q_{1}+q_{2}f\right\} \frac{ds}{ds^{\ast }}%
N+q_{2}^{\prime }\frac{ds}{ds^{\ast }}B
\end{equation*}%
Since $\widetilde{\gamma }^{\prime \prime }(s)$ has no components on $T$ and 
$B$, $q_{1}^{\prime }=0$ and $q_{2}^{\prime }=0.$Hence, $q_{1}$ and $q_{2}$
are constants. Thus, we prove the following theorem.

\begin{thm}
\label{th22}Let $\gamma $ and $\widetilde{\gamma }$ be the regular curves in 
$E^{3}.$If $\widetilde{\gamma }$ is the Bertrand mate curve of $\gamma $, $g$
is a constant.
\end{thm}

From the theorems \ref{th3} and \ref{th22}, we can conclude that there
exists the following relation between $g$ and $\widetilde{g}$ for a Bertrand
pair $\left( \gamma ,\widetilde{\gamma }\right) $ in $E^{3}$.%
\begin{equation*}
\epsilon g+\widetilde{g}=0.
\end{equation*}

\begin{cor}
Let $\gamma $ and $\widetilde{\gamma }$ be regular curves in $E^{3}.$ The
pair of $\left( \gamma ,\widetilde{\gamma }\right) $ is a Bertrand pair if
and only if $g$ and $\widetilde{g}$\ are constants.
\end{cor}

The tangent indicatrix of the Bertrand mate curve is%
\begin{equation*}
\widetilde{\gamma }_{t}=\frac{-1}{\sqrt{1+g^{2}}}\left\{ T-gB\right\} .
\end{equation*}

\begin{thm}
If the Frenet frame of the tangent indicatrix $\widetilde{\gamma }_{t}=%
\widetilde{T}$ of a Bertrand mate curve is $\{\widetilde{T}_{t},\widetilde{N}%
_{t},\widetilde{B}_{t}\}$, we have Frenet formula:%
\begin{equation*}
\widetilde{T}_{t}^{\prime }=\widetilde{\kappa }_{t}\widetilde{N}_{t},\text{
\ \ }\widetilde{N}_{t}^{\prime }=-\widetilde{\kappa }_{t}\widetilde{T}_{t}+%
\widetilde{\tau }_{t}\widetilde{B}_{t},\text{ \ \ }\widetilde{B}_{t}^{\prime
}=-\widetilde{\tau }_{t}\widetilde{N}_{t}
\end{equation*}%
where%
\begin{equation}
\widetilde{T}_{t}=-N\text{ \ , \ }\widetilde{N}_{t}=\frac{1}{\sqrt{1+f^{2}}}%
\left\{ T-fB\right\} \text{ \ , \ }\widetilde{B}_{t}=\frac{1}{\sqrt{1+f^{2}}}%
\left\{ fT+B\right\}  \label{16}
\end{equation}%
and%
\begin{equation}
s_{t}^{\ast }=-\dint \frac{\kappa \left( g-f\right) ^{2}}{f\left(
1+g^{2}\right) }ds^{\ast },\text{ \ }\widetilde{\kappa }_{t}=\frac{\sqrt{%
1+f^{2}}\sqrt{1+g^{2}}}{\left( f-g\right) },\text{ \ }\widetilde{\tau }_{t}=%
\frac{\kappa ^{\prime }\sqrt{1+g^{2}}}{\kappa ^{2}\left( 1+f^{2}\right) }
\label{17}
\end{equation}%
and $s_{t}^{\ast }$ is a natural representation of the tangent indicatrix of
the curve $\widetilde{\gamma }$. $\widetilde{\kappa }_{t}$ and $\widetilde{%
\tau }_{t}$ are the curvature and torsion of $\widetilde{\gamma }_{t}$,
respectively. The geodesic curvature of the principal image of the principal
normal indicatrix of $\widetilde{\gamma }_{t}$ is%
\begin{equation}
\widetilde{\Gamma }_{t}=.\frac{\kappa ^{3}\left( f-g\right) ^{2}\left(
1+f^{2}\right) ^{3/2}\left\{ \left( \kappa \kappa ^{\prime \prime }-3\kappa
^{\prime 2}\right) \left( 1+f^{2}\right) -3f\kappa ^{\prime 2}\left(
g-f\right) \right\} }{\sqrt{1+g^{2}}\left( \kappa ^{4}\left( 1+f^{2}\right)
^{3}+\kappa ^{\prime 2}\left( f-g\right) ^{2}\right) ^{3/2}}\frac{ds}{%
ds_{t}^{\ast }}.  \label{18}
\end{equation}
\end{thm}

\begin{thm}
\label{th25}Let $\left( \gamma ,\widetilde{\gamma }\right) $ be a
non-helical and non-planar Bertrand pair parameterized by the arclengths $s$
and $s^{\ast }.$Thus, the followings are true.

i. Bertrand mate curve is a slant helix if and only if the tangent
indicatrix of the Bertrand mate curve is a spherical helix.

ii. Bertrand curve is a slant helix if and only if the tangent indicatrix of
the Bertrand mate curve is a spherical helix.

iii. Bertrand curve is a slant helix if and only if the tangent indicatrix
of the Bertrand curve is a spherical helix.
\end{thm}

The arclengths $s_{t}^{\ast }$ and $s$ satisfy the equation%
\begin{equation*}
s_{t}^{\ast }=\dint \frac{\kappa \left( f-g\right) }{\sqrt{1+g^{2}}}ds.
\end{equation*}

We can re-state theorem \ref{th8} for Bertrand mate curve to characterize
its spherical image of tangent indicatrix as a helix without using $%
"\thicksim "$.

The principal normal indicatrix of the Bertrand mate curve is%
\begin{equation*}
\widetilde{\gamma }_{n}=\epsilon N.
\end{equation*}

\begin{thm}
If the Frenet frame of the principal normal indicatrix $\widetilde{\gamma }%
_{n}=\widetilde{N}$ of a non-helical and non-planar Bertrand mate curve is $%
\{\widetilde{T}_{n},\widetilde{N}_{n},\widetilde{B}_{n}\}$, we have Frenet
formula:%
\begin{equation*}
\widetilde{T}_{n}^{\prime }=\widetilde{\kappa }_{n}\widetilde{N}_{n},\text{
\ \ }\widetilde{N}_{n}^{\prime }=-\widetilde{\kappa }_{n}\widetilde{T}_{n}+%
\widetilde{\tau }_{n}\widetilde{B}_{n},\text{ \ \ }\widetilde{B}_{n}^{\prime
}=-\widetilde{\tau }_{n}\widetilde{N}_{n}
\end{equation*}%
where%
\begin{eqnarray}
\widetilde{T}_{n} &=&\frac{-\epsilon }{\sqrt{1+f^{2}}}\left\{ T-fB\right\} 
\notag \\
\widetilde{N}_{n} &=&\frac{\epsilon }{\sigma \sqrt{1+f^{2}}}\left\{ f\kappa
^{\prime }(g-f)T-\kappa ^{2}\left( 1+f^{2}\right) ^{2}N+\kappa ^{\prime
}(g-f)B\right\}  \label{19} \\
\widetilde{B}_{n} &=&\frac{1}{\sigma }\left\{ f\kappa ^{2}\left(
1+f^{2}\right) T+\kappa ^{\prime }(g-f)N+\kappa ^{2}\left( 1+f^{2}\right)
B\right\}  \notag
\end{eqnarray}%
where%
\begin{equation*}
\sigma =\sqrt{\kappa ^{\prime 2}\left( g-f\right) ^{2}+\kappa ^{4}\left(
1+f^{2}\right) ^{3}}
\end{equation*}%
and%
\begin{gather*}
s_{n}^{\ast }=\dint \frac{\kappa \left( g-f\right) \sqrt{1+f^{2}}}{f\sqrt{%
1+g^{2}}}ds^{\ast },\text{ }\ \ \widetilde{\kappa }_{n}=\frac{\sigma }{%
\kappa ^{2}\left( 1+f^{2}\right) ^{3/2}} \\
\widetilde{\tau }_{n}=\frac{-\epsilon (g-f)}{\sigma ^{2}}\left\{ 3f\kappa
^{\prime 2}\left( g-f\right) +\left( 3\kappa ^{\prime 2}-\kappa \kappa
^{\prime \prime }\right) \left( 1+f^{2}\right) \right\}
\end{gather*}%
and $s_{n}^{\ast }$ is a natural representation of the principal normal
indicatrix of the curve $\widetilde{\gamma }$. $\widetilde{\kappa }_{n}$ and 
$\widetilde{\tau }_{n}$ are the curvature and torsion of $\widetilde{\gamma }%
_{n}$, respectively.
\end{thm}

The arclengths $s_{n}^{\ast }$ and $s$ satisfy equation%
\begin{equation*}
s_{n}^{\ast }=\dint \kappa \sqrt{1+f^{2}}ds.
\end{equation*}%
We can re-state theorem \ref{th11} for Bertrand mate curve to characterize
its spherical image of principal indicatrix as a planar curve without using $%
"\thicksim "$.

The binormal indicatrix of the Bertrand mate curve is%
\begin{equation*}
\widetilde{\gamma }_{b}=\frac{-\epsilon }{\sqrt{1+g^{2}}}\left\{
gT+B\right\} .
\end{equation*}

If the Frenet frame of the binormal indicatrix $\widetilde{\gamma }_{b}=%
\widetilde{B}$ of a Bertrand mate curve is $\{\widetilde{T}_{b},\widetilde{N}%
_{b},\widetilde{B}_{b}\}$, we have Frenet formula:%
\begin{equation*}
\widetilde{T}_{b}^{\prime }=\widetilde{\kappa }_{b}\widetilde{N}_{b},\text{
\ \ }\widetilde{N}_{b}^{\prime }=-\widetilde{\kappa }_{b}\widetilde{T}_{b}+%
\widetilde{\tau }_{b}\widetilde{B}_{b},\text{ \ \ }\widetilde{B}_{b}^{\prime
}=-\widetilde{\tau }_{b}\widetilde{N}_{b}
\end{equation*}%
where%
\begin{equation}
\widetilde{T}_{b}=\epsilon N\text{, \ }\widetilde{N}_{b}=\frac{-\epsilon }{%
\sqrt{1+f^{2}}}\left\{ T-fB\right\} \text{, \ }\widetilde{B}_{b}=\frac{1}{%
\sqrt{1+f^{2}}}\left\{ fT+B\right\}  \label{20}
\end{equation}%
and%
\begin{equation}
s_{b}^{\ast }=-\dint \frac{\kappa \left( f-g\right) ^{2}}{f\left(
1+g^{2}\right) }ds^{\ast },\text{ \ \ }\widetilde{\kappa }_{b}=\frac{\sqrt{%
1+g^{2}}\sqrt{1+f^{2}}}{\left( f-g\right) },\text{ \ \ }\widetilde{\tau }%
_{b}=\frac{-\epsilon \kappa ^{\prime }\sqrt{1+g^{2}}}{\kappa ^{2}\left(
1+f^{2}\right) }  \label{21}
\end{equation}%
and $s_{b}^{\ast }$ is a natural representation of the binormal indicatrix
of the curve $\widetilde{\gamma }$ and it is equal to the total torsion of
the curve $\widetilde{\gamma }$. $\widetilde{\kappa }_{b}$ and $\widetilde{%
\tau }_{b}$ are the curvature and torsion of $\widetilde{\gamma }_{b}$,
respectively. The geodesic curvature of the principal image of the principal
normal indicatrix of $\widetilde{\gamma }_{b}$ is%
\begin{equation}
\widetilde{\Gamma }_{b}=\frac{\kappa ^{3}\left( f-g\right) ^{2}\left(
1+f^{2}\right) ^{3/2}\left\{ \left( \kappa \kappa ^{\prime \prime }-3\kappa
^{\prime 2}\right) \left( 1+f^{2}\right) -3f\kappa ^{\prime 2}\left(
g-f\right) \right\} }{\sqrt{1+g^{2}}\left( \kappa ^{4}\left( 1+f^{2}\right)
^{3}+\kappa ^{\prime 2}\left( f-g\right) ^{2}\right) ^{3/2}}\frac{ds}{%
ds_{b}^{\ast }}.  \label{22}
\end{equation}%
From equations (\ref{17}), (\ref{18}), (\ref{21}) and (\ref{22}), we
conclude that 
\begin{equation*}
\frac{\widetilde{\tau }_{t}}{\widetilde{\kappa }_{t}}=\frac{\widetilde{\tau }%
_{b}}{\widetilde{\kappa }_{b}}
\end{equation*}%
and $\widetilde{\Gamma }_{t}=\widetilde{\Gamma }_{b}.$

\begin{cor}
\label{cr33}The arclengths $s_{b}^{\ast }$ and $s$ satisfy the following
equation.%
\begin{equation*}
s_{b}^{\ast }=\dint \frac{\kappa \left( f-g\right) }{\sqrt{1+g^{2}}}ds
\end{equation*}
\end{cor}

Since $\lambda $\ is a constant, from (\ref{14}), we can write 
\begin{equation*}
\tau ^{\prime }=\lambda \kappa \kappa ^{\prime }(g-f)
\end{equation*}%
and since $q_{2}$ is a constant, we can write%
\begin{equation*}
\tau ^{\prime }=c_{1}\kappa ^{\prime }\sqrt{1+g^{2}}.
\end{equation*}%
Thus, we obtain%
\begin{equation*}
\frac{\kappa (g-f)}{\sqrt{1+g^{2}}}=\frac{c_{1}}{\lambda }
\end{equation*}%
and we can also say that%
\begin{equation*}
\frac{\kappa (g-f)}{\sqrt{1+g^{2}}}
\end{equation*}%
is a constant. From corollary \ref{cr33}, we obtain%
\begin{equation*}
\frac{ds_{b}^{\ast }}{ds}=\frac{c_{1}}{\lambda }
\end{equation*}%
and

\begin{equation*}
s_{b}^{\ast }=\frac{c_{1}}{\lambda }s+c_{2}
\end{equation*}%
where $c_{1}$ and $c_{2}$ are constants, too.

\begin{thm}
\label{teo33}Let $\left( \gamma ,\widetilde{\gamma }\right) $ be a
non-helical and non-planar Bertrand pairs parameterized by arclengths $s$
and $s^{\ast }.$Thus, the Bertrand curve is a slant helix if and only if the
binormal indicatrix of the Bertrand mate curve is a spherical helix.
\end{thm}

The theorem \ref{teo33} can easily prove as similar to theorem \ref{teo15}.
Furthermore, we can re-state theorem \ref{th17} for Bertrand mate curve to
characterize its spherical image of the binormal indicatrix as a spherical
helix without using $"\thicksim "$ and we can re-state corollary \ref{cr18}
for Bertrand mate curves.

We can give the following corollary by using (\ref{16}), (\ref{19}), and (%
\ref{20}).

\begin{cor}
The equations between the Frenet vectors of the spherical indicatrices of a
Bertrand mate curve in Euclidean 3-space are%
\begin{equation*}
\widetilde{T}_{t}=-\epsilon \widetilde{T}_{b},\text{ }\widetilde{N}%
_{t}=-\epsilon \widetilde{T}_{n}=-\epsilon \widetilde{N}_{b}\text{\ and }%
\widetilde{B}_{t}=\widetilde{B}_{b}.
\end{equation*}
\end{cor}

Both the spherical indicatrices of Bertrand curve and Bertrand mate curve
were examined in detail and it was seen that none of them didn't create any
specific curve pairs as a mutual.

\bigskip

\textbf{Acknowledgement}

We would like to thank the referee for valuable suggestions which improved
the paper considerably.


\begin{thebibliography}{9}
\bibitem{1} F. Wang and H. Liu, Mannheim Partner Curve in 3-Euclidean Space,
Mathematics in Practice and Theory, 37 (2007) 141-143.

\bibitem{2} H. Liu and F. Wang, Mannheim Partner Curve in 3-Space, Journal
of Geometry, 88 2008 120-126.

\bibitem{3} H. Matsuda and S. Yorozu, Notes on Bertrand Curves, Yokohama
Math. J., 50 (2003) 41-58.

\bibitem{4} L. Kula and Y. Yayli, On slant helix and its spherical
indicatrix, Applied Mathematics and Computation, 169 (2005) 600-607.

\bibitem{5} M.P. Carmo, Differential Geometry of Curves and Surfaces,
Prentice-Hall, New Jersey 1976.

\bibitem{6} S. Izumiya and N. Takeuchi, New Special Curves and Developable
Surfaces, Turk. J. Math., 28 (2004) 153-164.
\end{thebibliography}
\end{document}